\documentstyle[12pt,amscd,amssymb,epsf]{amsart} 
\newcommand{\p}{p}\newcommand{\q}{q}\newcommand{\sym}{S} 
\newcommand{\nse}{\kern-3pt\ns$=$}\newcommand{\qd}{\hfill\qed\medbreak} 
\newcommand{\To}{{\Bbb T}} 
\newcommand{\TT}{\widetilde T}\newcommand{\Z}{\mbox{$\Bbb Z$}} 
\newcommand{\C}{\mbox{$\Bbb C$}}\newcommand{\R}{\mbox{$\Bbb R$}} 
\newcommand{\noi}{\noindent}\newcommand{\qq}{\qquad\qquad} 
\newcommand{\ext}{\raise1pt\hbox{$\ts\bigwedge$}} 
\newcommand{\ds}{\displaystyle}\newcommand{\ts}{\textstyle} 
\newcommand{\rf}[1]{(\ref{#1})}\renewcommand{\c}[1]{\cite{#1}} 
\newcommand{\chii}{\raise2pt\hbox{$\chi$}}\newcommand{\Si}{\Sigma} 
\newcommand{\vep}{\varepsilon}\newcommand{\La}{\Lambda} 
\newcommand{\al}{\alpha}\newcommand{\bt}{\beta}\newcommand{\g}{\gamma} 
\newcommand{\F}{\mbox{$\cal F$}}\newcommand{\up}{\upsilon} 
\newcommand{\M}{\mbox{${\cal M}$}}
\newcommand{\Mg}{\mbox{${\cal M}\kern-2pt_g$}}
\newcommand{\Ng}{\mbox{${\cal N}\kern-2pt_g$}}
\newcommand{\U}{\mbox{${\cal U}$}}\newcommand{\V}{V\kern-1pt} 
\renewcommand{\d}{\hbox{\ns d}}\newcommand{\ddt}{\frac{\d}{\d t}}
\newcommand{\G}{\mbox{$\cal G$}}\newcommand{\I}{\mbox{$\cal I$}} 
\newcommand{\cir}{\raise1.6pt\hbox{\footnotesize$\circ$}} 
\newcommand{\tr}{\mbox{\bf tr}}\newcommand{\Vm}{\hbox{Vm}} 
\newcommand{\Res}[2]{\hbox{\ns Res}\kern-16pt\lower5pt
\hbox{\footnotesize$_{#1}$}\kern2pt\left[#2\right]}
\newcommand{\qk}{quaternion-K\"ahler}\renewcommand{\,}{\kern1pt} 
\newcommand{\ch}{\mbox{\ns ch}}\newcommand{\che}{\mbox{\ns ch}\,} 
\newcommand{\td}{\hbox{td}}  
\newcommand{\ns}{\normalshape}\newcommand{\ul}{\underline} 
\newcommand{\proof}{\noi{\it Proof.}\ }\newcommand{\op}{\oplus} 
\newcommand{\ot}{\otimes}  
\renewcommand{\ds}{\displaystyle}\renewcommand{\ts}{\textstyle} 
0 
\newtheorem{theo}{Theorem}[section]\newcommand{\A}{\mbox{$\widehat A$}} 
\newtheorem{lemma}{Lemma}[section] 
\newtheorem{corol}{Corollary}[section]\newcommand{\Ol}{\mbox{${\cal O}$}} 
\newtheorem{prop}{Proposition}[section]\def\frac#1#2{{#1\over#2}} 
\def\be#1\ee{\begin{equation}#1\end{equation}}
\begin{document} 

\title[Intersection Numbers and Symmetries of a Verlinde Formula] 
{Intersection Numbers on Moduli Spaces and Symmetries of a Verlinde Formula} 
\author{R.~Herrera\hspace{.5in} S.~Salamon } \address{Mathematical Institute 
\\24--29 St Giles'\\Oxford, OX1$\>$3LB}\thanks{$^*$The first author is 
supported by a scholarship from DGAPA, National University of Mexico} 
\date{25 December 1996}\maketitle\vspace{3pt}\bigbreak

\noi{\bf Abstract.} We investigate the geometry and topology of a
standard moduli space of stable bundles on a Riemann surface, and use
a generalization of the Verlinde formula to derive results on
intersection pairings.

\section[1]{Introduction}

Let $\Mg$ denote the smooth moduli space of stable holomorphic rank 2
vector bundles with fixed determinant of odd degree over a Riemann
surface $\Si_g$ of genus $g$. The space $\Mg$ has the structure of a
complex $(3g\!- \!3)$-dimensional K\"ahler manifold whose
anticanonical bundle is the square of an ample line bundle $L$
\cite{Ses,R}. The dimensions $h^0(\Mg,\Ol(L^j))$ are independent of
the complex structure on $\Si_g$ and were predicted in \c{V}. In this
paper, we highlight additional calculations that arise from the
Desale-Ramanan description \cite{DR} of $\Mg$ for the case in which
$\Si_g$ is hyperelliptic. Our approach is based on the proof of the
Verlinde formula by Szenes \c{Sz}, and grew out of an attempt to
generalize the related twistor geometry studied by the second author
in \c{parma}.

Universal cohomology classes $\al,\bt,\g$ were defined on $\Mg$ by Newstead
\c{N2} and used to compute the Chern character of the holomorphic tangent
bundle $T=T^{1,0}\Mg$. The latter can be expressed in terms of a tautological
rank $g-1$ bundle $Q$ with the help of the Adams operator $\psi^2$, and we show
in \S2 that the appearance of this algebraic gadget leads to quick proofs of
both the equation $\bt^g=0$ and the recurrence relation for the Chern classes
of $Q$. Underlying this theory is the fact that $\bt$ coincides with the
pullback of the canonical \qk\ class on a real Grassmannian, providing a link
with the index theory in \c{LS}. This approach has a number of simplifying
features, though in other ways \S2 is a supplement to the papers of Baranovsky
\cite{Ba} and Siebert-Tian \c{ST}. These authors, together with Zagier \c{Z}
and King-Newstead \c{KN}, have provided a variety of methods for determining
the cohomology ring of $\Mg$.

In \S3 we further exploit Adams operators and the fundamental role played by
$Q$ by computing the holomorphic Euler characteristics
\[\V_{g-1}(\p,\q)=\chii (\Mg,\Ol(\psi^{\p-\q}Q\ot L^{\q-1}))\] for all
$\p,\q\in\Z$, thereby extending the Verlinde formula (corresponding to $\p=\q$)
into a 2-dimensional array. The symmetries of the title are those enjoyed by
the integers $\V_{g-1}(\p,\q)$ in the $(\p,\q)$-plane. The main purpose of \S4
is to show that these symmetries can be used to recover the intersection
numbers $\big<\al^m\bt^n\g^p,[\Mg]\big>$, using calculations similar to those
of Thaddeus \c{Th} and Donaldson \c{D}. This formulation in turn provides an
alternative route to the Bertram-Szenes proof \c{BS} of the `untwisted'
Verlinde formula for the moduli space of semistable rank 2 bundles with fixed
determinant of even degree.

We end up encoding the topology of $\Mg$ into equations involving Chern
characters that are particularly easy to write down and remember. For example,
it is known that $\chii(\Mg,\Ol(T^*))=g-1$, and it follows from the
Riemann-Roch theorem that $\big<\che(\TT)\td(T),[\Mg]\big>=0$ where
$\TT=T^*-g+1$. We prove a stronger vanishing theorem, namely that
$\che(\TT)e^\alpha$ evaluates to zero when paired with any power of $\bt$. The
virtual bundle $\TT$ is a natural one to consider since it has virtual rank
$2g-2$ and $c_i(\TT)=0$ for $i>2g-2$ by \c{Gie}. Although our methods are very
special to the rank 2 case, it may be that analogues of these results hold on
moduli spaces of higher rank bundles.

\medbreak\section[1]{Tangent relations} 

Let $\Si_g$ be a hyperelliptic curve of genus $g$, admitting a double-covering
$\Si_g\to\C {\Bbb P}^1$ with distinct branching points $\omega_1,\ldots,
\omega_{2g+2}$. Desale and Ramanan proved in \c{DR} that the manifold $\Mg$
defined in \S1 can then be realized as the variety of $(g-1)$-dimensional
subspaces of $\C^{2g+2}$ isotropic with respect to the two quadratic forms
\be\sum_{i=1}^{2g+2}y_i^2,\qquad\sum_{i=1}^{2g+2}\omega_i y_i^2.\label{qf}\ee
One therefore obtains a holomorphic embedding of $\Mg$ into the complex
homogeneous space \[\F_g=\frac{SO(2g+2)} {U(g-1)\times SO(4)}\] parametrizing
$(g-1)$-dimensional subspaces which are isotropic with respect to the first
quadratic form \c{Sz}.

Let $Q,\>W$ denote the duals of the tautological complex vector bundles over 
$\F_g$ with fibres $\C^{g-1},\>\C^4$ and structure groups $U(g-1),\>SO(4)$ 
respectively. (The notation $Q$ is consistent with \c{ST}, and from now on we 
shall often drop the subscript $_g$ since the genus will generally be fixed in 
our discussion.) The decomposition of the standard representation of $SO(2g+2)$ 
on $\C^{2g+2}$ under $U(g-1)\times SO(4)$ provides the equation \be Q^*\op Q\op 
W=2g+2,\label{q*} \ee where the right-hand side denotes a trivial vector bundle 
of rank $2g+2$. The second form in \rf{qf} now determines a non-degenerate 
section $s$ of the symmetrized tensor product $\sym^2Q$, and the zero set of 
$s$ coincides with $\M$. 

The holomorphic tangent bundle $T^{1,0}\F$ of $\F$ is determined by the summand
$\frak m$ in the Lie algebra splitting \[{\frak s\frak o}(2g+2)_c\cong({\frak
u}(g-1)\op{\frak s\frak o}(4))_c\op{\frak m}\op\overline{\frak m},\] where $_c$
denotes complexification. Given that ${\frak s\frak o}(2g+2)_c\cong
\ext^2\C^{2g+2}$, it follows from \rf{q*} that we may choose the orientation so
that \[T^{1,0}\F\cong\ext^2Q\>\op\>(Q\ot W).\] This equation is well known in
the context of twistor spaces, since $\F$ is a 3-symmetric twistor space that
fibres over the oriented real Grassmannian \be\G_g={\Bbb G}{\ns r}_4(\R^{2g+2})
=\frac{SO(2g+2)}{SO(2g-2)\times SO(4)}\label{Gr}\ee for $g\ge3$ \c{Br,1164}.
The term $\ext^2Q$ is simply the holomorphic tangent bundle to the Hermitian
symmetric fibres $SO(2g-2)/U(g-1)$, and its complement $Q\ot W$ corresponds to
the holomorphic horizontal bundle that plays an important role in the theory of
harmonic maps \cite{BR}.

With the above choices, the normal bundle of $\M$ in $\F$ is isomorphic to 
$\sym^2Q$, and \[T^{1,0}\F\vert_{\cal M}\cong T^{1,0}\M\>\op\>\sym^2Q 
\vert_{\cal M}.\] In the notation of K-theory, we may write \[T=T^{1,0}\M= 
\ext^2Q\>+\>QW\>-\>\sym^2Q,\] where from now on we are using the same symbols 
to denote bundles pulled back to $\M$. Writing $\psi^2=\sym^2-\ext^2$, we have 

\begin{lemma} $T=QW-\psi^2Q$.\end{lemma} 

\noi The operator $\psi^2$ is one of the series of Adams operators, defined by
the formula \[\sum_{p\ge0}(\psi^p\!E)t^p=r-t\ddt\log\La_{-t}E,\] where $E\in
K(\M)$ has virtual rank $r$ and $\La_tE=\sum_{i\ge0}(\ext^iE)t^i$ \c{FL}. Each
$\psi^p$ is a ring homomorphism in K-theory, and is characterized by the
property that \be\ch_k(\psi^p\!E)=p^k\ch_k(E),\label{Adam} \ee where $\ch_k$
denotes the term of dimension $2k$ in the Chern character. We shall use the
operators $\psi^p$ with $p\ge3$ in \S3. \medbreak

Cohomology classes \be\al\in H^2(\M,{\Bbb Z}),\quad\bt\in H^4(\M,{\Bbb Z}),
\quad \g\in H^6(\M,{\Bbb Z})\label{uc}\ee were introduced by Newstead \c{N,AB}.
They are obtained from the K\"unneth components of the characteristic class
$c_2({\Bbb V})$, where $\Bbb V$ is a universal $SO(3)$ bundle over $\M$, and
generate the ring $H^*_I(\M)$ of cohomology classes of $\M$ invariant by the
action of the mapping class group on $H^3(\M)$. By expressing $T$ in terms of a
push-forward of $\Bbb V$ one obtains the following result, which is effectively
the definition of \rf{uc} for our purposes:

\smallbreak\begin{theo}{\nse\cite[Theorem~2]{N2}} \[\che(T)=3g-3+2\al+ 
\sum_{k\ge2}\frac{s_k}{k!},\quad \hbox{\ns where}\quad\left\{\begin{array}{rcl} 
s_{2k-1}&=&2\al\bt^{k-1}-8(k-1)\g\bt^{k-2},\\[3pt]s_{2k}&=&2(g-1)\bt^k. 
\end{array}\right.\]\end{theo} \smallskip 

As an application of Lemma 2.1 and \rf{q*}, we see that the complexification of 
the real tangent bundle of $\M$ is isomorphic to \begin{eqnarray}T+T^*&=& 
(Q^*+Q)W\>-\>\psi^2(Q^*+Q),\nonumber\\&=&(2g+2-W)W\>-\>(2g+2-\psi^2W)\nonumber 
\\&=&(2g+2)(W-1)-W^2+\psi^2W.\label{TM}\end{eqnarray} The Chern character of 
this may be read off and then compared with Theorem~2.1 (for this purpose it 
helps not to replace $-W^2+\psi^2W$ by the equivalent $-2\ext^2W$). An easy 
calculation gives \be\;\che(T+T^*)=6g\!-\!6+2(g\!-\!1)p_1(W)+\hbox{$\frac16$} 
\big[(g\!-\!1)p_1(W)^2-2(g\!+\!5)p_2(W)\big]+\cdots\label{pp}\ee where the 
Pontrjagin classes are defined by regarding $W$ as an $SO(4)$-bundle. Now 
\rf{pp} must equal twice the sum of the even terms of $\che(T)$, so $p_1(W)= 
\bt$, $p_2(W)=0$ and \be\che(W)=2+ e^{\sqrt\bt}+ e^{-\sqrt\bt}\label{cosh}\ee 
on the moduli space $\M$. 

Since $Q^*+Q=2g+2-W$ is a genuine complex vector bundle of rank $2g-2$ with 
total Chern class \[c(W)^{-1}=\sum_{k=0}^\infty\bt^k\] on $\M$, we get 
Conjecture~(a) of \c{N2}: 

\begin{theo}\label{bg=0} $\bt^g=0$.\end{theo} 

\noi This was first proved by Kirwan \c{K}, who established the completeness of 
the Mumford relations on $H^*_I(\M)$. It also follows from results in 
\c{JW,Wi}, and a different proof was given by Weitsman \c{W} in the more 
general setting of moduli spaces of flat connections over a Riemann surface 
with marked points.\medbreak 

{}From Lemma~2.1 and Theorem~2.1 one may readily compute the Chern
character of $Q$ in terms of the classes \rf{uc}. From this point of
view the definition of $\M$ as a submanifold of $\F$ could not be
simpler, as the terms $\ext^2Q$ and $\sym^2Q$ `miraculously' combine
into a form specifically adapted for computing characters. The result
is

\begin{theo} \[\che(Q)=g-1+\al+\sum_{k\ge2}\frac{s'_k}{k!},\quad\hbox{\ns 
where}\quad\left\{\begin{array}{rcl}s'_{2k-1}&=&\al\bt^{k-1}+2\g\bt^{k-2}, 
\\[3pt]s'_{2k}&=&-\bt^k.\end{array}\right.\]\end{theo} 
            
\proof Let $s_k,\,s'_k$ denote the components of $\che(T),\,\che(Q)$,
respectively, in dimension $2k$. Using Lemma~2.1, \rf{Adam} and
\rf{cosh},
\[3g-3+\sum_{k\ge1}{s_k\over k!}=2\Big(g-1+\sum_{k\ge1}{s'_k\over k!}\Big) 
\Big(2+\sum_{k\ge1}{\bt^k\over(2k)!}\Big)-\Big(g-1+\sum_{k\ge1}{2^k
s'_k\over k!}\Big).\] The result now follows from Theorem~2.1 by
induction on $k$.\qd
                                                
An analogue of the last equation can be found in \c{Ba}, though the authors 
were led to it by the paper of Siebert and Tian \cite{ST}, who give an 
equivalent expression for $\che(Q)$. Theorem~2.3 leads very quickly to their 
recurrence relation for the Chern classes of $Q$. Using a standard trick \c{Z}, 
the generating function for the Chern classes of $Q$ is recaptured by the 
formula \begin{eqnarray}c(t)&=&\exp\Big[\sum_{k\ge1}{(-1)^{k-1}s_k t^k\over k} 
\Big]\nonumber\\[2pt]&=&\exp\Big[\al t+\sum_{n\ge2}(\al\bt^{n-1}+2\g\bt^{n-2}) 
\frac{t^{2n-1}}{2n-1}+\sum_{n\ge1}\bt^n\frac{t^{2n}}{2n}\Big].\label{ct} 
\end{eqnarray} The relation \cite[Proposition~25]{ST}, namely \[(1-\bt 
t^2)c'(t)=(\al+\bt t +2\g t^2)c(t)\] follows immediately, whence 
               
\begin{corol} The Chern classes of the rank $g-1$ bundle $Q$ on $\Mg$ satisfy 
\[(k+1)c_{k+1}=\al c_k+k\bt c_{k-1}+2\g c_{k-2}.\]\end{corol}\bigbreak 

Zagier has shown that the $c_k$ coincide with the K\"unneth components of the
Chern classes in $H^*(\Mg)\ot H^*(J_g)$ of a higher rank bundle used to define
the Mumford relations ($J_g$ is the Jacobian of $\Si_g$). The identities in
$\al,\bt,\g$ arising from the equations $c_k=0$ for $k=g,g+1,g+2$ then provide
a minimal set of relations to completely determine the cohomology ring
$H_I^*(\M)$ \cite{Z,Ba,KN,ST}. We have set out to show that these equations
follow quite directly from Newstead's own computations, and it is worth
pointing out that Corollary~2.1 is analogous, but simpler, to the recurrence
relation for the Chern classes of $T$ given at the end of \cite{N2}.\medbreak

The fact that the Pontrjagin ring of $\M=\Mg$ is generated by $p_1(\M)=2(g-1) 
\bt$ can also be related to the geometry of the real Grassmannian \rf{Gr}. For 
$Q+Q^*$ and $W$ are (complexifications of) the pullbacks of real vector bundles 
$\hat U$, $\hat W$ over $\G=\G_g$, and the real tangent bundle of $\G$ is 
isomorphic to $\hat U\ot\hat W$. The choice of an orientation of $W$ gives the 
manifold $\G$ a \qk\ structure. The latter is characterized by a certain 
non-degenerate closed 4-form $\Omega$ that arises from the curvature of a 
locally-defined quaternionic line bundle $H$, and represents the integral 
cohomology class $4u\in H^4(\G,\Z)$, where $u=-c_2(H)$ \c{LS}. 

\begin{prop} $\bt$ is the pullback of the class $4u$ by means of the embedding 
cum projection $\M\hookrightarrow\F\to\G$.\end{prop} 

\proof From above, $\bt$ is the pullback to $\M$ of $\hat\bt=p_1(\hat W)$. A 
calculation from \c{parma} shows that \[ p_2(\hat W)=(\hat\bt-4u)^2\in H^8( 
\G,\Z).\] Assuming that $g\ge3$, $b_4(\M)=2$ and so $\hat\bt-4u$ must pull back 
to $a\al^2+b\bt$ on $\M$ for some $a,b\in\Z$; from \rf{cosh}, 
$(a\al^2+b\bt)^2=0$. However, the remarks after Corollary~2.1 imply that there 
are no non-trivial relations involving $\al^4,\,\al^2\bt,\,\bt^2$ in $H^8(\M)$ 
except that \[0=c_4(Q)-\al\,c_3(Q)=\a^4+2\al^2\bt-3\bt^2\] in genus 3. (There 
are actually four distinct \qk\ structures on $\G_3=SO(8)/(SO(4)\times SO(4))$, 
and Proposition~2.1 holds for only two of them.) It follows that in all cases 
$\bt=4u$ in $H^4(\M)$.\qd\vspace{.1in}

Any \qk\ manifold $M$ of dimension $4m$ is the base space of a complex manifold
$Z$ (the more usual `twistor space') fibred by rational curves. The positive
integer \be\up(M)=\big<(4u)^m,[M]\big>=\frac1{2(m+1)^{2m+1}}
\,c_1(Z)\sp{2m+1}\ee determines the `quaternionic volume' of $M$, and can be
expressed in terms of dimensions of representations of the isometry group,
using techniques from \c{LS}. For $M=\G_g$ we have $m=2g-2$, $Z\cong
SO(2g+2)/(SO(2g-2)\times U(2))$, and one can prove that
\[\up(\G_g)=\frac2g{4g-3\choose2g-1};\] by comparison $\up(\Bbb{HP}^{2g-2})=
4^{2g-2}$. Theorem~2.2 is in contrast to the non-degenerate nature of the
4-form $\Omega$ over $\G$, and reflects the failure of $\M$ to map onto a
quaternionic subvariety of $\G$.\medbreak 

\section[1]{Character calculations}\vspace{3pt} 
                
In this section, we set $h=g-1$ and consider the holomorphic Euler
characteristics \be\V_h(\p,\q)=\chii({\cal M}_{h +1},\Ol(\psi^{\p-\q}Q\ot
L^{\q-1})).\label{ab}\ee Given that $c_1(\M)=2\al$ and the canonical bundle of
$\M$ is isomorphic to $L^{-2}$, Serre duality implies that
\be\V_h(\p,\q)=(-1)^h\V_h(-\p,-\q),
\label{Sd} \ee with the convention that $\psi^{-\p}Q=\psi^\p Q^*$. Following 
\c{Sz}, we set $w=x+x^{-1}-2$ and \be F(w,\p)={(x^\p-x^{-\p})(x-x^{-1})
\over x+x^{-1}-2} =\sum_{h\ge 0}\left[4{\p+h\choose2h+1}+{\p+h-1\choose2h-1}
\right]w^h,\label{fwk}\ee in order to define 
\[{F(w,\p)\over F(w,\q)^2}=\sum_{k\ge0} G_k(\p,\q)w^k.\]
The next theorem expresses \rf{ab} in terms of this generating function.

\begin{theo}\label{extension} Let $h\ge2$. Then $\V_h(\p,0)=0$ for all $\p$, 
and
\[\V_h(\p,\q)=4(-4\q)^h\left(\p(h+1)\,G_h(\q,\q)-\q\,G_h(\p,\q)\right),
\quad\p>0.\] In particular, $\V_h(\p,\q)+\V_h(-\p,\q)=0$ for all
$\p,\q\in\Z$.\end{theo}\vspace{4pt}

The resulting symmetries of $\V_h(\p,\q)$ are illustrated schematically in
Figure~1.\vspace{5pt}

\begin{corol}\label{p=cq} Let $c$ be an integer such that $0\le c\le2+(g-2)/
\q$. Then \be\V_h(c\,\q,\q)=c\sum_{j=1}^\q(-1)^{j+1}(h+1-(-1)^{j(c+1)})
\left({\q\over\sin^2(j\pi/2\q)}\right)^h.\label{caa}\ee\end{corol}
\vspace{7pt}
                                  
Setting $\p=\q$ in \rf{ab} gives \[\V_h(\q,\q)=h\dim H^0(\M,\Ol(L^{\q-1})),\]
since the higher cohomology spaces are zero by Kodaira vanishing and Serre
duality. When $c=1$, the right-hand side of \rf{caa} does indeed reduce to the
Verlinde formula for the dimension of the space of sections of $L^{\q-1}$. This
was first deduced from fusion rules \c{V}; a direct proof was given by Szenes
\c{Sz}, though many other generalizations now exist
\c{NRs,JW,Th2,BL,Fa}. Moreover, the right-hand side of \rf{fwk} is essentially
the generating function for the reciprocal of the Verlinde series. Note also
that \be {1\over h}\V_h(1,1)=\chii(\M,\Ol)=\big<\td(\M), [\M]\big>=1
\label{Todd}\ee is the Todd genus of $\M$. Finally, when $h=1$ we have $Q\cong
L$ and it follows that $\V_1(\p,\q)=\V_1(\p,\p)$ for all $\p,\q$.

\vspace{-15pt}\begin{center}\leavevmode \epsfxsize260pt\epsfbox[10 0 350 
350]{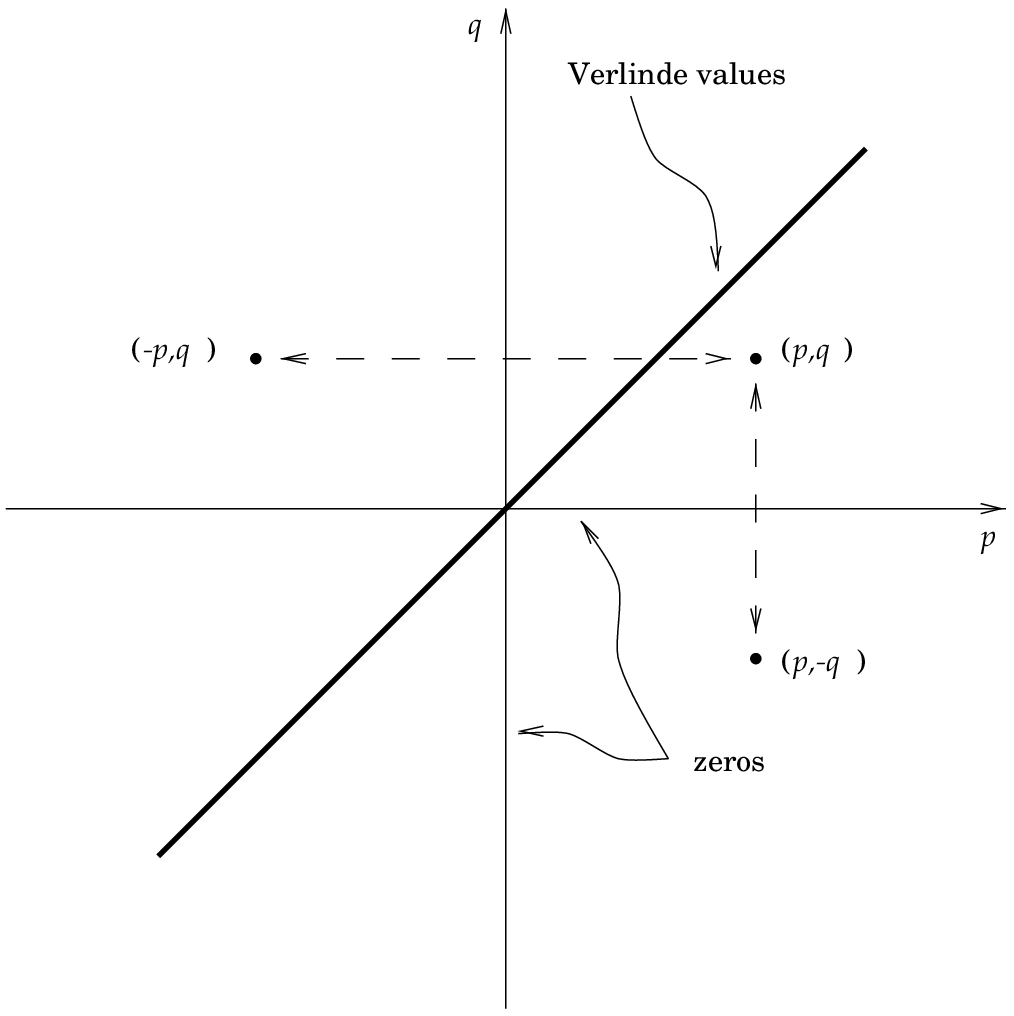}\end{center}\centerline{Figure 1}\bigbreak

\vspace{.2in}\noi{\it Proof of Theorem~3.1.} We follow closely Szenes' proof of
the Verlinde formula in \c{Sz}, and his use of the

\begin{lemma}{\nse\cite[Proposition~4.3]{FL}} Let $E$ be a vector bundle of 
rank $n$ over a smooth projective variety $X$, and let $i\colon M\to X$ be the 
zero locus of a non-degenerate section of $E$. Then $\chii(M,\Ol(i^*U))=\chii 
(X,\Ol(U\ot\La_{-1}E^*))$ for any vector bundle $U$ on 
$X$.\end{lemma}\smallbreak 

On a homogeneous space, holomorphic Euler characteristics can be computed by 
means of the Atiyah-Bott fixed point formula. Let $G$ be a reductive Lie group, 
$P$ a parabolic subgroup of $G$, and $F=G/P$ the corresponding flag manifold. A 
representation $R$ of $P$ determines both a holomorphic vector bundle $\ul 
R=G\times_P R$ over $F$ and a virtual $G$-module \[\I(R)=\sum_i(-1)^i H^i(F, 
\Ol(\ul R)).\] Let $\To$ be a common maximal torus of $P$ and $G$, let 
$W_G,\,W_P$ be the Weyl groups, and $W_r$ the relative Weyl group. The 
character of the $G$-module $\I_R$ is given by \be\tr(\I_R)= \sum_{w\in W_r}w 
\cdot{\tr(R)\over\tr(\La_{-1}A^*)},\label{trace}\ee where $A$ is the $P$-module 
associated to the holomorphic tangent bundle $T^{1,0}\F=\ul A$. The right-hand 
side of \rf{trace} is a function on $\To$, and evaluation at the identity 
element yields \[\tr(\I_R)\vert_e=\chii(F,\Ol (\ul R)).\]\medbreak 

Returning to the problem at hand, let $H$ denote the subgroup $U(h)\times
SO(4)$ of $G=SO(2g+2)$, where $h=g-1$. If $B$ denotes the fundamental
representation of $U(h)$, then the vector bundles $\ul B^*$ and $\det\ul B^*$
over $\F=SO(2g+2)/H$ pull back to $Q$ and $L$ respectively over $\M$. Lemma~3.1
therefore implies that
\[\chii(\M,\Ol(\psi^{\p-\q}Q\ot L^{\q-1}))=\chii(\F,\Ol(\ul R)),\]
where \[R=R^{\p,\q}=\psi^{\p-\q}B^*\ot(\det B^*)^{\q-1}\ot\La_{-1}
(\sym^2B).\]

Now we proceed to calculate (\ref{trace}). Let $x_1,\ldots,x_{g+1}$ be the 
characters of the maximal torus of $SO(2g+2)$ corresponding to the polarisation 
$\{y_{2j-1}+iy_{2j}:1\le j\le g+1\}$ of $\C^{2g+2}$. The character of the 
fundamental $SO(2g+2)$-module is $\sum_{j=1}^{h+2}(x_j+x_j^{-1})$, and that of 
the fundamental $U(h)$-module $\sum_{j=1}^h x_j^{-1}$. Thus, 

\[\tr(R^{\p,\q}) 
=\prod_{i\le h}x_i^{\q-1}\prod_{1\le j\le k\le h}\Big(1-{1\over x_jx_k} 
\Big)\sum_{\ell=1}^h x_\ell^{\p-\q},\] and from \c{Sz}, \[\tr(\La_{-1} 
A^*)=\prod_{1\le i<j\le h}\bigg(1-{1\over x_ix_j}\bigg)\prod_{1\le k\le h\atop 
\vep=1,2}\left(1-{1\over x_{h+\vep}x_k}\right)\left(1-{x_{h+\vep}\over x_k} 
\right),\] 

\noi so it suffices to prove that \be\V_h(\p,\q)=\lim_{\{x_i\to 1\}} \sum_{w\in
W_r}w\cdot{\tr(R^{\p,\q})\over\tr(\La_{-1}A^*)}.\label{lim}\ee\smallbreak

We have that 

\[{\tr(R^{\p,\q})\over\tr(\La_{-1}A^*)}=\left(\prod_{i\le h}
{x_i^\q(x_i-x_i^{-1})\over(x_i+x_i^{-1}-x_{h+1}-x_{h+1}^{-1})(x_i+x_i^{-1}-
x_{h+2}-x_{h+2}^{-1})}\right)\,\sum_{j\le h}x_j^{\p-\q}.\] 

\noi To perform the summation in \rf{lim} as in \c{Sz} we first recall the form
of the relative Weyl group $W_r$ of $W_{SO(2h+4)}$ with respect to $W_{U(h)}$
and $W_{SO(4)}$. It has $2^h{h+2\choose2}$ elements, and

\[W_r=W^{\ns signs}\rtimes W^{perms}\] 

\noi where $W^{signs}$ consists of all the substitutions $x_i\mapsto x_i^{-1}$
of an {\it even} number of variables modulo $\{x_{h+1}\mapsto x_{h
+1}^{-1},x_{h+2}\mapsto x_{h+2}^{-1}\}$, and $W^{perms}$ consists of all the
permutations of two variables.  Adding up first with respect to $W^{signs}$ we
get

{\small\[\left(\prod_{i\le h}{(x_i- x_i^{-1})\over(x_i+x_i^{-1}-x_{h
+1}-x_{h+1}^{-1})(x_i+x_i^{-1}-x_{h+2}-x_{h+ 2}^{-1})}\right)\,\sum_{j\le h}
(x_j^\p-x_j^{-\p})\,\prod_{l\ne j}(x_j^\q- x_j^{-\q}),\]}

\noi setting $x_{h+2} \mapsto 1 $ and then adding up with respect to
$W^{perms}$ gives a contribution

{\small\[\sum_{j=1}^{h+1}\left[\left( \prod_{i\ne
j}{x_i-x_i^{-1}\over(x_i+x_i^{-1}-2)(x_i+x_i^{-1}-x_j-x_j^{-1})}
\right)\,\sum_{k\ne j}\left((x_k^\p-x_k^{-\p})\,\prod_{\ell\ne j\atop\ell\ne k}
(x_\ell^\q-x_\ell^{-\q})\right)\right].\]}\medbreak

Substituting $w_i=x_i+x_i^{-1}-2$, and rearranging the terms converts the above
summation into 

{\small\be\left(\sum_{j=1}^{h+1}F(w_j,\p)\, \prod_{i\ne
j}F(w_i,\q)\right)\,\sum_{j=1}^{h+1}{(-1)^{j-1}\Vm(w_1,\ldots, \widehat
w_j,\ldots,w_{h+1})\over F(w_j,\q)\,\Vm(w_1,\ldots,w_{h+1})}\qq\qq\label{exp2}
\ee\vspace{-4pt}\[\qq\qq\qq-\left(\prod_{i=1}^{h+1}F(w_i,\q)\right)\, 
\sum_{j=1}^{h+1}{(-1)^{j-1}F(w_j,\p)\Vm(w_1,\ldots,\widehat w_j,\ldots,w_{h+1})
\over F(w_j,\q)^2\,\Vm(w_1,\ldots,w_{h+1})},\]} 

\noi where $\Vm$ denotes the Vandermonde determinant. Since
$\lim\limits_{x\to1}F(w,\q)=4\q$, the first factors in both summands converge
as the $x_i$ tend to $1$. By \cite[Lemma~5.3]{Sz} the first summand tends to
$(-4\q)^h4\p(h+1)G_h(\q,\q)$, and the second to $(-4\q)^{h+1}G_h(\p,\q)$.\qd

\bigbreak\noi{\it Proof of Corollary 3.1.} The hypothesis on $c$ implies that
the meromorphic form \[{F(w,c\q)dw\over F(w,\q)^2w^{h+1}}\] over $\Bbb{CP}^1$
has no poles at $0$ and $\infty$. The result is then a consequence the residue
theorem and \cite[Lemma~5.3]{Sz}.\qd \vfil\eject

\medbreak\section[1]{Further relations} 

The identity \be\V(\p,0)=0,\qquad p\in\Z,\label{i}\ee of Theorem~3.1 is also an
easy consequence of Theorems~2.2, 2.3 and the Hirzebruch-Riemann-Roch
theorem. For the latter implies that
\[\V(\p,0)=\big<\che(\psi^\p Q\ot L^*)\td(\M),[\M]\big>=\big<\che(\psi^\p Q)
\A(\M),[\M]\big>,\] and the $\A$ class \be\A(\M)=\left({{1\over2}
\sqrt\bt\over\sinh{1\over2}\sqrt\bt}\right)^{2g-2} \kern-10pt,\label{Ahat}\ee
readily computed from \rf{TM} and \rf{cosh}, is a polynomial in $\bt$.

The identity \rf{Ahat} was used by Thaddeus to show that the Verlinde formula 
(Corollary~3.1 with $c=1$) determines the intersection form on $\M$. In 
\cite[Equation~(30)]{Th}, he gives the intersection numbers \be\big<\al^m\bt^n 
\g^p,[\Mg]\big> =(-1)^{p-g}{g!\,m!\over(g-p)!q!}2^{2g-2-p}(2^q-2)B_q,\label{IN} 
\ee where $m+2n+3p=3g-3$, $q=m+p-g+1$, and $B_q$ is the $q$th Bernoulli number 
(equal to $q!$ times the coefficient of $x^q$ in $x/(e^x-1)$). Another key 
point in the argument is \cite[Proposition~(26)]{Th}, namely that $\g$ is 
Poincar\'e dual to $2g$ copies of ${\cal M}_{g-1}$, so that \be\big<\al^m\bt^n 
\g^p,[\Mg]\big>=2g\big<\al^m\bt^n\g^{p-1},[{\cal M}_{g-1}]\big>,\qquad m+2n+3p= 
3g-3.\label{26}\ee\bigbreak 

In this section, we shall show that the intersection numbers \rf{IN} are in
fact determined by \rf{i} and the identities \begin{eqnarray} && \qquad
V_h(0,\p)=0,\label{ii}\\&&V_h(\p,\p)+\V(-\p,\p)=0,\quad \p\in\Z\label{iii}
\end{eqnarray} that follow from Theorem~3.1, or rather its proof. 
Following closely the notation of Donaldson 
\cite[\S5]{D}, we set \[I_{k}^{(g)}={1\over(g-1+2k)!}\big<\al^{g-1+2k}
\bt^{g-1-k},[\M]\big>,\] and also \[I^{(g)}(t)=\sum_{k=0}^{g-1}I_k^{(g)}
t^{2k}.\]

\begin{theo} $\ds I^{(g)}(t)=(-4)^{g-1}{t\over\sinh t}$.\end{theo} 

\proof Interpreting \rf{i} as a polynomial identity in $\p$, using the 
Hirzebruch-Riemann-Roch theorem and \rf{26}, yields the equation \be 
t^2\ddt\left({I^{(g)}(t)\,\sinh t\over t}\right) =g(t-\sinh t)(I^{(g)}(t)+ 
4I^{(g-1)}(t))\label{aa}\ee modulo $t^{2g}$. Similarly, \rf{iii} gives 
\[t^2\ddt\left({I^{(g)}(t)\,\sinh\,2t\over t}\right)+t\,I^{(g)}(t) 
(1-\cosh\,2t)= g(2t-\sinh\,2t)(I^{(g)}(t)+4I^{(g-1)}(t)),\] which can
be simplified into \be 2t^2\ddt\left({I^{(g)}(t)\sinh t\over
t}\right)=g(2t- \sinh(2t))(I^{(g)}(t)+ 4I^{(g-1)}(t)).\label{cc}\ee
Both sides of \rf{aa} and \rf{cc} must now be identically zero, and
ignoring the modulo $t^{2g}$, $I^{(g)}(t)=C(g)t/\sinh t$ where $C(g)$
is a function of $g$ such that $C(g)+4C(g-1)=0$. But, using the
description \rf{qf} of ${\cal M}_2$ as the intersection of two
quadrics in $\Bbb{CP}^5$, we find that $C(2)=-\big<\al^3, [{\cal
M}_2]\big>=-4$. It also follows now that
$I^{(g)}(t)=(-4)^{g-1}K(t^2)$, where $K$ is the generating function of
the intersection pairings described in \c{D}.\qd\vspace{7pt}

Let $\Ng$ denote the moduli space of semistable rank 2 vector bundles with
fixed even degree determinant over a compact Riemann surface $\Si_g$ of genus
$g$, and let $L_0$ be the generator of $\hbox{\ns Pic}(\Ng)$. An easy
application of the above knowledge of $I^{(g)}(t)$ now yields

\begin{theo}{\nse\cite[Theorem~1.1]{BS}} $\;\ds\dim H^0 (\Ng,L_0^{k-2})= 
\sum_{j=1}^{k-1}\left({k\over 1-\cos(2j\pi/k)}\right)^{g-1}$. \end{theo} 

\proof Let $\U$ be a universal rank 2 bundle over $\Mg\times\Si_g$, so
that if $m\in\M$ represents the bundle $E\to\Si$ then
$\U\vert_{\{m\}\times\Si}\cong E$; this may be chosen such that if
$U_x=\U\vert_{\M\times\{x\}}$ then $\det(U_x)\cong L$ over
$\M$. Bertram and Szenes \c{BS} use a Hecke correspondence to prove
that \[\dim H^0(\Ng,L_0^k)=\chii(\Mg,\sym^kU_x).\] Note that $c(U_x)=
1+\al+\frac14(\al^2-\bt)$, so that $c(U_x\ot L^{-1/2})=1-\bt/4$, and
\[\che(\sym^k (U_x\ot L^{-1/2}))={\sinh((k+1)\sqrt{\bt}/2)\over
\sinh(\sqrt{\bt}/2)}\] \c{JW} Then we need to evaluate
\[\left<{\sinh((k+1)\sqrt{\bt}/2)\over\sqrt{\bt}/2}
e^{(k+2)\alpha/2}\left({\sqrt{\bt}/2\over\sinh(\sqrt{\bt}/2)}\right)^{2g-1}
\kern-8pt,\;[\M]\right>.\] Using Theorem~4.1, this is readily seen to equal
\begin{eqnarray}&&\kern-6pt\Res{t=0}{\left({k+2\over-2\sinh^2(t/(k+1))}
\right)^{g-1}\!\left({k+1\over k+2}\right){\sinh
t\over\sinh((k+2)t/(k+1))\sinh(t/(k+1))}dt}
\nonumber\\&&=\sum_{j=1}^{k+1}\left({k+2\over 1-\cos(2\pi j/(k+2))}
\right)^{g-1}\kern-4pt\sqrt{-1}\kern5pt\Res{\!u=\pi
j}{{\sin((k+1)u/(k+2))\over \sin
u\sin(u/(k+2))}du}\nonumber\\&&=\sum_{j=1}^{k+1}\left({k+2\over
1-\cos(2\pi
j/(k+2))}\right)^{g-1}.\qquad\qed\nonumber\end{eqnarray}\vspace{7pt}

In the last part of this paper, we shall combine separate uses of the
Adams operators, namely Lemma~2.1 and the identities \rf{ii},\rf{iii}
to express the intersection form of the smooth moduli space $\M=\Mg$
in an alternative form. The bundle $Q$ was defined geometrically only
in the hyperelliptic case, and to de-emphasise its role at this stage
we consider in addition \[\TT=T^*-g+1\in K(\M).\] Note that this has
virtual rank $2g-2$ and vanishing higher Chern classes $c_i(\TT)=0$
for $i>2g-2$ by Gieseker's theorem \cite{Gie,Z}.

Let us say that a cohomology class $\delta\in H^*(\M)$ is {\it saturated} if 
$\big<\delta\,\bt^j,[\M]\big>=0$ for all $j\ge0$. Thus, any polynomial in $\bt$ 
is itself saturated. With this terminology, 

\begin{prop} $\che(Q^*\ot L)$ and $\che(\TT\ot L)$ are both saturated, 
i.e.\par\noi {\ns(a)} $\big<\che(Q^*)e^\alpha\bt^j,[\M]\big>=0,\quad j\ge0$; 
\par\noi {\ns(b)} $\big<\che(\TT)e^\alpha\bt^j,[\M]\big>=0,\quad j\ge0$. 
\end{prop} 

\proof Since $\psi^\p$ is a ring homomorphism in K-theory, \rf{ii} implies 
\[0=\big<\che(\psi^\p Q^*\ot L^{\p-1})\td(\M),[\M]\big>= 
\big<\che(\psi^\p(Q^*\ot L))\A(\M),[\M]\big>.\] Equation (a) follows from the 
fact that the identity above is true for {\it all} $\p\in{\Bbb Z}$. Now 
consider the decomposition \[(T^*-g+1)\ot L=Q^*\ot L\ot W-(\psi^2Q^*+ g-1)\ot 
L.\] \rf{cosh} and part (a) imply that $\che(Q^*\ot L\ot W)$ is saturated. It 
therefore suffices to show that $\che((\psi^2Q^*+g-1)\ot L)$ is saturated, but 
given that \[\psi^\p((\psi^2Q+g-1)\ot L^*)=\psi^{2\p}Q\ot L^{-\p}\>+\>(g-1) 
L^{-\p},\] this follows from \rf{iii}.\qd 

The equations of Proposition~4.1 for $j\ge g-1$ follow immediately from 
Theorems~2.2 and 2.3, though taking $j=0,\ldots,g-2$ gives independent 
relations. For example, expanding (a) shows that 
\[2k\,I^{(g)}_k+(g+2k)\sum_{i=1}^k {I^{(g)}_{k-i}\over(2i+1)!} \] is a linear 
combination of pairings $\big<\al^m\bt^n\g,[\M]\big>$ for each $k$ with $1\le 
k\le g-1$. The intersection numbers {\ns\rf{IN}} can then be determined by the 
equations of Proposition~4.1(a) (with Theorem~2.3), \rf{Todd} and \rf{26} by 
induction on the genus $g$. 

In view of the Riemann-Roch equation \[\chii(\M,T^*)-g+1=\big<\che(\TT) 
e^\alpha\hat A(\M),[\M]\big>\] that follows from \rf{Todd}, Proposition~4.1(b) 
correctly predicts the coefficient \[\chii(\M,T^*)=g-1\label{norm}\] of $t$ in 
the polynomial $\chii_t=(1-t)^{g-1}(1+t)^{2g-2}$ computed in \cite{NR}. This 
observation suggests that there should exist more direct proofs of 
Proposition~4.1. 

\medbreak\newcommand{\bi}{\bibitem} 

\enddocument